\documentclass[12pt]{article}
\usepackage{amsmath}
\usepackage{amsfonts}
\usepackage{amssymb}
\usepackage{amsthm}
\usepackage{hyperref}

\usepackage{graphicx}

\begin{document}

\thispagestyle{empty}
\centerline{}
\centerline{\Large\bf All solutions of the Diophantine equation }
\centerline{\Large\bf $a^6+b^6=c^6+d^6+e^6+f^6+g^6$ }
\centerline{\Large\bf for $a,b,c,d,e,f,g<250000$ }
\centerline{\Large\bf found with a distributed Boinc project }

\bigskip

\begin{center}

{\large \sl Robert Gerbicz}\\
{\large \sl Jean-Charles Meyrignac}\\
{\large \sl Uwe Beckert}\\*[5mm]
\bigskip
August, 2011
\end{center}

\bigskip

\begin{center}
{\bf Abstract}\\
\end{center}
The above equation is also called as Euler(6,2,5) system.
By computational aspect these systems are very interesting. And we can also apply these methods to other Diophantine equations. We give a brief history of these systems and how we searched for these big solutions on Boinc. Our two Boinc projects ran from April of $2010$ to July of $2011$.
\bigskip
\begin{minipage}{12.8cm}

\end{minipage}

\bigskip

\section{History}
In $1753$, the famous mathematician Leonhard Euler conjectured that it would take {\bf at least} $k$ $k$th powers to sum to an $k$th power. This is known as the {\bf Euler's Sum of Powers Conjecture}\cite{ESPC1}\cite{ESPC2}.\\
The first example for $k>3$ turned up only after more than $150$ years, when R. Norrie finally found in $1911$,\\
$30^4+120^4+272^4+315^4=353^4$.\\
\\
Surprisingly, this conjecture turned out to be {\bf false} and is one instance where Euler was wrong.
In $1966$, L. J. Lander and T. R. Parkin found with a computer search that $27^5+84^5+110^5+133^5=144^5$ \cite{LanderParkin}.\\
\\
In $1998$, Randy L. Ekl defined {\bf Euler's Extended Conjecture}\cite{Ekl}. This can be simply stated as:\\
\\
The diophantine equation $(k,m,n)$,\\
\begin{center}
$x_1^k+x_2^k+\cdots +x_m^k=y_1^k+y_2^k+\cdots +y_n^k$
\end{center}
has no (positive) integer solution for $k>m+n$, other than the trivial case when all $x_i=y_i$ (or its permutation).
\\
\\
Here is a list of the most important discovered solutions, in chronological order:\\
$(4,2,2)$ $158^4+59^4=134^4+133^4$ (L. Euler, $1722$)\\
$(6,3,3)$ $23^6+15^6+10^6=22^6+19^6+3^6$ (K. Subba Rao, $1934$)\\
$(5,1,4)$ $144^5=133^5+110^5+84^5+27^5$ (L. J. Lander and T. R. Parkin, $1966$)\\
$(4,1,3)$ $20615673^4=18796760^4+15365639^4+2682440^4$ (Noam Elkies, $1986$)\\
$(4,1,3)$ $422481^4=414560^4+217519^4+95800^4$ (Roger Frye, $1988$)\\
$(5,2,3)$ $14132^5+220^5=14068^5+6237^5+5027^5$ (Bob Scher and Ed Seidl, $1997$)\\
$(8,3,5)$ $966^8+539^8+81^8=954^8+725^8+481^8+310^8+158^8$ (Scott I. Chase, $2000$)\\
$(5,1,4)$ $85359^5=85282^5+28969^5+3183^5+55^5$ (Jim Frye, $08/27/2004$)\\
$(8,4,4)$ $3113^8+2012^8+1953^8+861^8=2823^8+2767^8+2557^8+1128^8$ (Nuutti Kuosa, $11/09/2006$)\\
\\
The current most wanted solutions are $(6,1,5)$ and $(6,2,4)$.
\\
In $1999$, Jean-Charles Meyrignac started a distributed project \cite{EULER} to enumerate all solutions to $(6,2,5)$, in the hope to discover a solution to $(6,1,5)$ or $(6,2,4)$, if one of the terms is zero.\\
\\
The algorithm was designed by Giovanni Resta, and required only a few megabytes of memory.\\
The smallest known solution was discovered independently by Edward Brisse and Giovanni Resta in $1999$:\\
$(6,2,5)$ $1117^6+770^6=1092^6+861^6+602^6+212^6+84^6$.
\\
After $11$ years of computation (and almost $310$ years of total computer time), all solutions with the biggest term less than $88000$ were discovered.
\\
In $2008$, Robert Gerbicz wrote a much faster program, requiring a larger memory set.

\section{Computation}
D. J. Bernstein in $1999$ wrote a great article about how to enumerate solutions of\\
 $p(a)+q(b)=r(c)+s(d)$, in this article he applied his method to various Euler systems. I've used a different way, in spite of this it is worth to read Bernstein's classic article. \cite{Bernstein}

To start with from Fermat's little theorem it is easy to see that $x^{p-1}\pmod p$ is one, if $x$ is relative prime to $p$, and zero if $x$ is divisible by $p$ (here $p$ is prime). For $p=7$ this gives that $a^6+b^6\equiv u_7 \pmod p$, where $u_7$ is the number of terms of $\{a,b\}$ which are not divisible by $7$. Therefore  $c^6+d^6+e^6+f^6+g^6\equiv u_7\pmod 7$, hence on the right side there are also $u_7$ terms which are relative prime to $7$. So at least three terms of these are divisible by $7$  (because $u_7\le 2$), let's call these $e,f,g$. Note that $u_7=0$ can't be, otherwise all terms are divisible by $7$, and that wouldn't be a primitive solution; so at least one of $a,b$ and at least one of $c,d$ are relative prime to $7$.\\
Hence $(a^6+b^6)-(c^6+d^6)=e^6+f^6+g^6$, here $e,f,g$ are divisible by $7$, it gives that $(a^6+b^6)-(c^6+d^6)$ is divisible by $7^6=117649$. In modular tricks this gives us the largest speedup.\\
\\
In addition to this we can get more modular speedup. This time by Euler-Fermat theorem: $x^{\varphi(m)}\equiv 1\pmod m$ if $gcd(x,m)=1$, apply this for $m=9$: $x^6\equiv 1\pmod 9$ if $gcd(x,9)=1$, and trivially: $x^6\equiv 0\pmod 9$ if $gcd(x,9)>1$. Again by a similar proof as above we get that at least three terms of the right side are divisible by $3$, and this gives that at least one of the $e,f,g$ terms is divisible by three.\\
It isn't hard to see that $x^6\pmod 8$ is one if $x$ is odd and zero if $x$ is even. We can obtain again that at least three terms on the right side is even, so at least one of the $e,f,g$ terms is even.\\
\\
Moreover there are "deeper" modular tricks, I will show one. Note that\\ $(a^6+b^6)-(c^6+d^6)=7^6*((\frac e7)^6+(\frac f7)^6+(\frac g7)^6)$, we learned that $x^6\equiv \{0,1\}\pmod 7$, from this:\\ $(a^6+b^6)-(c^6+d^6)\equiv 7^6*\{0,1,2,3\}\pmod {7^7}$.\\

Anyway the question remains how we searched for solutions? Forget all the above.
There is an $O(N^4*\log N)$ solution: generate all $e^6+f^6+g^6$ triplet in a table, then for all $a,b,c,d$ quadruplet check by binary search if $(a^6+b^6)-(c^6+d^6)$ is in the table or not, if it is in the table then print out the solution. Is it manageable? Yes, if $N$ is small, but in our search $N=250000$ and to build up this table it takes $O(N^3)$ memory, and that's too large.\\
\\
Better method: choose a $p$ prime number near $N^{1.5}$ and consider the equation modulo $p$. Fix $0\le r_p<p$ and generate all triplets, where $e^6+f^6+g^6\equiv r_p\pmod p$. For each $0\le w<p$ remainder do the following: first generate all $(a,b)$ pairs for that $a^6+b^6\equiv w\pmod p$. We can do it in $O(N)$ time: for each $0\le a<N$ we know that $b^6\equiv w-a^6\pmod p$, and if $p$ is in form $3k+2$ then this has got at most two solutions (modulo $p$), so by a precomputed table we can get them in $O(1)$ time. Similarly generate all solutions of $c^6+d^6\equiv w-r_p\pmod p$ in $O(N)$ time. We can expect that the number of $(a,b)$ (and $(c,d)$) pairs is $O(N^2/p)=O(\sqrt N)$. Subsequently  by a double loop in $O(N*\log N)$ time we can check if the $(a^6+b^6)-(c^6+d^6)$ sum is in the stored $e^6+f^6+g^6$ table or not. How fast is it? The same as above, because for $w$ and $r_p$ we have $p$ possibilities, this gives $O(p^2*N*\log N)=O(N^4*\log N)$ time, but it uses only $O(N^3/p)=O(N^{1.5})$ memory to store the table for $e^6+f^6+g^6\equiv r_p\pmod p$. Note that for each $r_p$ value we can generate this table in $O(N^2)$ time by the precomputed table.\\

The consequence is that in memory we are good, but we lost the above modular tricks. The great news is that we can use all of them!\\
For $e,f,g$ triplets it is easy to see that if $e^6+f^6+g^6\equiv r_p\pmod p$ then $(\frac e7)^6+(\frac f7)^6+(\frac g7)^6\equiv \frac{r_p}{7^6}\pmod p$, and here $0\le \frac e7,\frac f7,\frac g7<\frac N7$ are integers. Furthermore at least one of the $e,f,g$ is even and at least one of them is divisible by three.
To use the factor of $7^6$: fix $0<v<7^6$ and generate all $(a,b)$ pairs for that $a^6+b^6\equiv v\pmod {7^6}$, and place this pair to the $T_s$ list, where $(a^6+b^6)\equiv s\pmod p$. Use these lists to get solutions of $(a^6+b^6)-(c^6+d^6)\equiv r_p\pmod p$: one solution from $T_s$ and one from $T_{s-r_p}$. I must also point out that we can use the ``deeper`` tricks, and in the code we avoid the binary searches, for this reason we get a speedup by $O(\log N)$. Using the big factor of $7^6$ modify also the size of $p$.\\

All in all the computation time to determine all solutions up to $N$ is $O(N^4)$ (with a very small multiplier), and it uses $O(N^{1.5})$ memory.

\section{Summary}
In practise our first Boinc search up to $N=117649$ took $70$ GHz years (this is measured on AMD Phenom processor) and used $128$MB of Ram. We found $32$ new primitive solutions and rediscovered all previously known solutions. The second bigger search up to $N=250000$ took about $740$ GHz years (in other words this means more than $2*10^{19}$ floating point operations) with an improved code and used at most $540$MB of Ram. Note that we computed each workunit at least twice to ensure that we don't miss a solution, so a single search would take about half of the above time. This time we found $196$ new solutions, thus up to $N=250000$ there are $377$ primitive solutions of the Euler$(6,2,5)$ system.

Thanks to those $4470$ users who computed at least one workunit!\\
Special thanks go to Uwe Beckert\cite{Rechenkraft} for providing and operating the Boinc\\
distributed volunteer project.

\section {Appendix} (primitive solutions up to $250000$ and the discoverer's name, \\the solutions are sorted in lexicographical order)

\begin{verbatim}
1. 1117^6+770^6=1092^6+861^6+602^6+212^6+84^6
   (Greg Childers)
2. 2041^6+691^6=1893^6+1468^6+1407^6+1302^6+1246^6
   (Greg Childers)
3. 2441^6+752^6=2184^6+2096^6+1484^6+1266^6+1239^6
   (Greg Childers)
4. 2827^6+151^6=2653^6+2296^6+1488^6+1281^6+390^6
   (Greg Childers)
5. 2959^6+2470^6=2954^6+2481^6+850^6+798^6+420^6
   (Greg Childers)
6. 6623^6+323^6=6615^6+2912^6+642^6+434^6+363^6
   (Marcin Lipinski)
7. 8905^6+347^6=8820^6+5489^6+2576^6+2499^6+534^6
   (tnerual)
8. 8969^6+5203^6=9023^6+2604^6+2520^6+2379^6+478^6
   (tnerual)
9. 9551^6+9451^6=10080^6+8589^6+4884^6+3976^6+3943^6
   (JC Meyrignac@work)
10. 9612^6+7271^6=9198^6+7446^6+6580^6+6279^6+5118^6
    (tnerual)
11. 9707^6+6277^6=9675^6+5796^6+5531^6+4536^6+3640^6
    (JC Meyrignac@work)
12. 12272^6+2459^6=11522^6+9144^6+8283^6+7434^6+1400^6
    (tnerual)
13. 13417^6+5933^6=11403^6+11378^6+10698^6+4641^6+70^6
    (Greg Childers)
14. 13903^6+5317^6=13788^6+8295^6+6026^6+1827^6+1232^6
    (Greg Childers)
15. 15149^6+6914^6=14520^6+10675^6+10444^6+6006^6+810^6
    (Greg Childers)
16. 15627^6+1865^6=14196^6+13608^6+5385^6+4032^6+565^6
    (Laurent Evrard)
17. 15905^6+10519^6=14679^6+13759^6+9534^6+3822^6+2482^6
    (tnerual)
18. 16103^6+9111^6=14868^6+13755^6+8484^6+6253^6+2004^6
    (Greg Childers)
19. 16159^6+10532^6=15442^6+11760^6+11016^6+10273^6+3234^6
    (Rich Brown)
20. 16160^6+1865^6=15792^6+11470^6+4984^6+4329^6+1008^6
    (Rich Brown)
21. 16315^6+8749^6=15246^6+11682^6+11375^6+11140^6+7119^6
    (Greg Childers)
22. 16867^6+14786^6=17790^6+9821^6+9786^6+2912^6+2502^6
    (Rich Brown)
23. 17737^6+16043^6=17646^6+14574^6+14261^6+3255^6+2590^6
    (Rich Brown)
24. 19168^6+1747^6=18984^6+10584^6+10088^6+8214^6+5299^6
    (tnerual)
25. 19869^6+10673^6=19866^6+9783^6+9324^6+2088^6+1925^6
    (Lars Dausch)
26. 20111^6+15059^6=18879^6+15680^6+14854^6+13812^6+2499^6
    (Mark Dodrill)
27. 20173^6+15459^6=20469^6+13118^6+11550^6+5355^6+1170^6
    (JC Meyrignac@home)
28. 20235^6+1690^6=19915^6+11706^6+11472^6+10542^6+5754^6
    (JC Meyrignac@work)
29. 20273^6+17527^6=20755^6+16182^6+8892^6+5348^6+4221^6
    (Laurent Evrard)
30. 23192^6+20259^6=22932^6+19170^6+16874^6+13629^6+7602^6
    (Lars Dausch)
31. 24781^6+125^6=23457^6+20046^6+6951^6+2968^6+214^6
    (Mark Dodrill)
32. 25153^6+6036^6=25116^6+11297^6+6798^6+6090^6+6006^6
    (user)
33. 25309^6+19441^6=23996^6+22302^6+11532^6+8897^6+2067^6
    (Rich Brown)
34. 25360^6+4407^6=25020^6+15372^6+13524^6+10542^6+5015^6
    (F. Yura)
35. 25951^6+4901^6=25750^6+15486^6+7203^6+3906^6+3031^6
    (Physics3)
36. 26143^6+21591^6=27132^6+16625^6+8871^6+5622^6+966^6
    (Rich Brown)
37. 26305^6+6344^6=25032^6+20940^6+9968^6+6303^6+1022^6
    (Physics10)
38. 26435^6+21913^6=27702^6+5887^6+3339^6+2568^6+1274^6
    (Rich Brown)
39. 26797^6+12368^6=24794^6+22764^6+11539^6+5148^6+1386^6
    (Rich Brown)
40. 27081^6+15688^6=25116^6+20762^6+20664^6+8322^6+6135^6
    (Mark Dodrill)
41. 28314^6+3211^6=27531^6+19006^6+17580^6+12096^6+4452^6
    (Lars Dausch)
42. 28537^6+2038^6=28112^6+16702^6+16602^6+12519^6+3360^6
    (JC Meyrignac@work)
43. 29040^6+22271^6=29804^6+16590^6+5748^6+4914^6+2625^6
    (Rich Brown)
44. 29384^6+17515^6=28518^6+22088^6+16128^6+9828^6+1355^6
    (Laurent Evrard)
45. 29875^6+28496^6=31584^6+25158^6+7536^6+2632^6+1901^6
    (Mark Dodrill)
46. 30038^6+4607^6=26726^6+26310^6+17650^6+12831^6+12474^6
    (Laurent Evrard)
47. 31291^6+848^6=31164^6+16818^6+1918^6+805^6+54^6
    (Mark Dodrill)
48. 31881^6+25204^6=30030^6+28644^6+14998^6+14484^6+4809^6
    (Andrew Taylor)
49. 32549^6+9575^6=28749^6+28518^6+20072^6+16506^6+10991^6
    (Laurent Evrard)
50. 32555^6+30271^6=34479^6+25116^6+16037^6+15078^6+2404^6
    (Mark Binfield)
51. 32831^6+13273^6=32612^6+19425^6+8967^6+8946^6+3946^6
    (Andrew Taylor)
52. 33566^6+14671^6=31090^6+28497^6+10752^6+4542^6+2800^6
    (Laurent Evrard)
53. 33835^6+4957^6=33097^6+22176^6+19978^6+12249^6+6552^6
    (Mark Dodrill)
54. 34016^6+16481^6=30174^6+30100^6+20223^6+11544^6+2324^6
    (Mark Dodrill)
55. 34359^6+22445^6=34374^6+20258^6+19467^6+1134^6+639^6
    (F. Yura)
56. 34459^6+10175^6=33754^6+21849^6+20972^6+11661^6+4872^6
    (peterzal)
57. 34649^6+23615^6=34392^6+24752^6+16331^6+6699^6+1614^6
    (Lars Dausch)
58. 35272^6+7769^6=34722^6+21952^6+18585^6+16506^6+7436^6
    (Laurent Evrard)
59. 38869^6+19799^6=36820^6+31689^6+12726^6+5578^6+4491^6
    (F. Yura)
60. 39217^6+10327^6=38886^6+22315^6+19670^6+7644^6+2709^6
    (Mark Dodrill)
61. 39791^6+4755^6=35742^6+32208^6+29481^6+21931^6+8316^6
    (Mark Dodrill)
62. 39902^6+26353^6=39126^6+25872^6+25137^6+24778^6+3010^6
    (Glenn Klakring)
63. 39955^6+34207^6=42114^6+19600^6+18057^6+7035^6+5152^6
    (michael allen)
64. 40261^6+31165^6=41489^6+20472^6+9642^6+8351^6+7224^6
    (Physics3)
65. 42013^6+27691^6=39725^6+31776^6+31254^6+19712^6+5733^6
    (Laurent Evrard)
66. 43387^6+3917^6=37401^6+37296^6+30594^6+27305^6+14000^6
    (h.d.)
67. 43691^6+12449^6=39039^6+38810^6+11676^6+5323^6+3234^6
    (vignobles and Mangraa)
68. 44081^6+22273^6=40187^6+36099^6+30394^6+25020^6+6072^6
    (Rich Brown)
69. 44911^6+16553^6=43680^6+31935^6+23527^6+19226^6+6720^6
    (F. Yura)
70. 45046^6+25775^6=40678^6+39676^6+22194^6+21168^6+13839^6
    (Andrew Taylor)
71. 45253^6+3664^6=44805^6+27762^6+18274^6+11592^6+7364^6
    (Mark Dodrill)
72. 45305^6+20873^6=42819^6+34118^6+31563^6+8300^6+7950^6
    (Mark Dodrill)
73. 45305^6+33654^6=42126^6+40647^6+10584^6+7154^6+4878^6
    (Mark Dodrill)
74. 46205^6+17570^6=41853^6+37128^6+34188^6+23482^6+10388^6
    (Laurent Evrard)
75. 47081^6+25308^6=47250^6+15977^6+13146^6+11898^6+10626^6
    (Laurent Evrard)
76. 47654^6+43619^6=48034^6+42402^6+26403^6+23422^6+252^6
    (Mark Dodrill)
77. 48149^6+9313^6=43337^6+39725^6+34944^6+20718^6+13704^6
    (h.d.)
78. 49085^6+39617^6=47942^6+41084^6+31038^6+13881^6+5733^6
    (Mark Dodrill)
79. 49999^6+44630^6=47906^6+47138^6+27810^6+13020^6+2289^6
    (IceFox)
80. 50264^6+38639^6=47658^6+44072^6+23982^6+23823^6+18214^6
    (Physics7)
81. 51179^6+16312^6=51051^6+22554^6+21238^6+19970^6+3582^6
    (Physics3)
82. 51293^6+16262^6=51296^6+14404^6+11418^6+9786^6+6531^6
    (Laurent Evrard)
83. 51299^6+34079^6=44954^6+41487^6+40019^6+36414^6+4110^6
    (Mark Dodrill)
84. 51761^6+35968^6=52689^6+10404^6+7444^6+4886^6+3864^6
    (Laurent Evrard)
85. 51928^6+21739^6=50694^6+37254^6+18837^6+16750^6+7822^6
    (Physics2)
86. 52733^6+25675^6=52297^6+32334^6+21966^6+20538^6+10343^6
    (Mark Dodrill)
87. 53573^6+32533^6=52122^6+35679^6+34454^6+31639^6+18324^6
    (IceFox)
88. 53603^6+34847^6=48385^6+44772^6+40474^6+24402^6+16323^6
    (Mark Dodrill)
89. 54093^6+22073^6=49126^6+42171^6+38241^6+36204^6+21822^6
    (Lars Dausch Desktop)
90. 55315^6+44707^6=56707^6+36684^6+30562^6+22386^6+7923^6
    (Lars Dausch Desktop)
91. 55623^6+8786^6=51996^6+46254^6+19758^6+11865^6+7252^6
    (toto)
92. 55623^6+35341^6=53589^6+42388^6+35490^6+20559^6+13608^6
    (toto)
93. 56303^6+30151^6=52734^6+41664^6+38267^6+37387^6+1260^6
    (Laurent Evrard)
94. 57299^6+11945^6=54258^6+41951^6+40523^6+6246^6+2646^6
    (S. I . Rodriguez)
95. 57509^6+47251^6=51485^6+49644^6+47544^6+35957^6+3744^6
    (Mark Dodrill)
96. 58013^6+44545^6=57876^6+42546^6+36632^6+10977^6+6307^6
    (Monarcho)
97. 58843^6+20179^6=52434^6+51569^6+31656^6+31409^6+16926^6
    (Laurent Lucas)
98. 59035^6+36959^6=58884^6+38318^6+16961^6+12933^6+8646^6
    (Bommer 2)
99. 59107^6+8635^6=53203^6+52038^6+21749^6+6678^6+4644^6
    (Lars Dausch Notebook)
100. 59995^6+6856^6=57974^6+41250^6+39372^6+15540^6+9499^6
     (Laurent Evrard)
101. 60323^6+4548^6=59108^6+36216^6+34944^6+33411^6+20274^6
     (tnerual)
102. 61815^6+15131^6=61041^6+36003^6+30534^6+29442^6+27524^6
     (Laurent Evrard)
103. 64729^6+4232^6=64218^6+38745^6+17866^6+8638^6+546^6
     (Laurent Evrard)
104. 66667^6+9013^6=63084^6+48699^6+46086^6+35035^6+3680^6
     (F. Yura)
105. 66999^6+26135^6=62342^6+56169^6+29484^6+10743^6+9114^6
     (desmao3)
106. 68963^6+8035^6=65163^6+54132^6+36428^6+35286^6+34139^6
     (Laurent Evrard)
107. 69311^6+21142^6=66906^6+47866^6+45426^6+27636^6+6529^6
     (Fritz)
108. 70873^6+51697^6=69289^6+52542^6+49164^6+9919^6+1596^6
     (Lars Dausch Games)
109. 70955^6+28561^6=63798^6+59037^6+50316^6+35959^6+2456^6
     (naoki)
110. 71593^6+14335^6=66675^6+58451^6+42588^6+31410^6+6202^6
     (bagleyd)
111. 71751^6+29341^6=63728^6+60291^6+52473^6+31482^6+24612^6
     (bagleyd)
112. 72458^6+39793^6=69708^6+54978^6+42420^6+28273^6+7856^6
     (Lars Dausch Games)
113. 72939^6+7228^6=71040^6+48888^6+44268^6+30261^6+21446^6
     (Laurent Evrard)
114. 72943^6+39643^6=66640^6+62850^6+40761^6+29834^6+1821^6
     (Lars Dausch Notebook)
115. 73772^6+69919^6=70266^6+64652^6+64449^6+48030^6+29822^6
     (Toutouf)
116. 74398^6+35411^6=71106^6+58422^6+34765^6+28884^6+23968^6
     (davidwillis90)
117. 74705^6+7075^6=65982^6+64302^6+50665^6+39291^6+16354^6
     (JC)
118. 75451^6+37005^6=73750^6+52038^6+42876^6+21105^6+9009^6
     (Laurent Evrard)
119. 77012^6+64091^6=79150^6+54096^6+43547^6+23310^6+15660^6
     (Cousin Caterpillar)
120. 77171^6+75180^6=85302^6+40478^6+35574^6+20979^6+6006^6
     (Laurent Evrard and JC)
121. 77797^6+70250^6=78680^6+68101^6+39858^6+28140^6+27162^6
     (bagleyd)
122. 77920^6+2753^6=71652^6+58968^6+54138^6+52639^6+4904^6
     (bagleyd)
123. 78322^6+1975^6=72212^6+58982^6+57015^6+48204^6+18660^6
     (Laurent Evrard)
124. 78943^6+39667^6=77676^6+48556^6+48531^6+22113^6+3346^6
     (bagleyd)
125. 79235^6+42305^6=73958^6+64236^6+51254^6+32193^6+19005^6
     (Lars Dausch Desktop)
126. 79612^6+4065^6=77826^6+55496^6+38367^6+11658^6+9132^6
     (Madcat)
127. 79621^6+35323^6=78990^6+47208^6+37401^6+15526^6+4577^6
     (Lars Dausch Notebook)
128. 79810^6+48939^6=77168^6+59598^6+49812^6+31152^6+10965^6
     (Jwb52z)
129. 79846^6+77017^6=75768^6+70924^6+68313^6+60640^6+3606^6
     (KAMCOBILL)
130. 80744^6+2749^6=77532^6+57850^6+53027^6+23940^6+9408^6
     (Laurent Evrard)
131. 80821^6+72010^6=83264^6+66024^6+34764^6+25963^6+6048^6
     (Laurent Evrard)
132. 80939^6+58025^6=76461^6+69356^6+41420^6+35217^6+32790^6
     (kinshuk)
133. 81757^6+76663^6=81972^6+68514^6+65051^6+51681^6+2114^6
     (Vesuvius)
134. 82067^6+60035^6=80526^6+65562^6+25025^6+4487^6+2010^6
     (Madcat)
135. 82822^6+4667^6=79758^6+55568^6+52965^6+48258^6+32536^6
     (bagleyd)
136. 82828^6+64207^6=84945^6+46652^6+43468^6+22596^6+20958^6
     (bagleyd)
137. 83057^6+41900^6=83202^6+32529^6+30072^6+17908^6+3916^6
     (JC)
138. 84521^6+44353^6=84312^6+48349^6+24285^6+13776^6+11648^6
     (JC)
139. 84673^6+19766^6=84252^6+40607^6+39742^6+36582^6+21468^6
     (steinrar)
140. 84685^6+40567^6=82515^6+61992^6+30212^6+22066^6+16695^6
     (LocalBusinessMen)
141. 84746^6+23003^6=73252^6+72891^6+63308^6+34626^6+11970^6
     (Laurent Evrard)
142. 86105^6+50615^6=83772^6+64071^6+45416^6+30345^6+12862^6
     (Laurent Evrard and JC)
143. 88063^6+64745^6=85776^6+62650^6+61945^6+53487^6+33432^6
     (JC)
144. 88153^6+80407^6=88864^6+75675^6+62370^6+27181^6+2436^6
     (1fast6)
145. 88352^6+43825^6=82256^6+65519^6+62286^6+57330^6+7434^6
     (Roald)
146. 88615^6+62283^6=86520^6+70392^6+33660^6+15351^6+10043^6
     (Roald)
147. 89857^6+28717^6=85501^6+71460^6+36540^6+28686^6+23569^6
     (Lars Dausch Games)
148. 90061^6+1562^6=80430^6+80038^6+10394^6+5421^6+4662^6
     (Pwrguru)
149. 90283^6+27359^6=89838^6+48937^6+36330^6+24166^6+13299^6
     (Lars Dausch Games)
150. 92711^6+47567^6=83027^6+80556^6+59802^6+14700^6+14029^6
     (Liuqyn and [AF>HFR>RR] Jim PROFIT)
151. 94082^6+46435^6=89166^6+75390^6+49746^6+35581^6+23206^6
     (Greeri and tr0ach)
152. 94196^6+66781^6=95517^6+54714^6+31584^6+5306^6+140^6
     (Liuqyn and proteino.de)
153. 95699^6+80957^6=91905^6+81991^6+71302^6+47754^6+618^6
     (Jeff17 and [AF>HFR>RR] julien76100)
154. 95713^6+63016^6=91080^6+79423^6+46074^6+9646^6+3402^6
     (Raetha and Liuqyn)
155. 96373^6+84914^6=95340^6+78741^6+72870^6+57550^6+29288^6
     (Iki and Jeff17)
156. 96977^6+1382^6=93583^6+71994^6+50820^6+39158^6+10308^6
     (Michael H.W. Weber and proteino.de)
157. 97057^6+30026^6=97062^6+26369^6+22698^6+10650^6+4354^6
     (MX-6 and Rebirther)
158. 97136^6+66217^6=94674^6+76300^6+42882^6+29673^6+13516^6
     (sdl* and Laurent Evrard)
159. 97159^6+64893^6=89222^6+83967^6+57141^6+54474^6+5490^6
     (Liuqyn and proteino.de)
160. 97243^6+32386^6=85211^6+80616^6+75738^6+29400^6+18368^6
     (proteino.de and Fire$torm.SETI.USA [BlackOps])
161. 97973^6+68227^6=87048^6+80031^6+69281^6+67494^6+65954^6
     (Laurent Evrard and Jeff17)
162. 98395^6+28723^6=94633^6+75712^6+33642^6+16089^6+594^6
     (Death and Jeff17)
163. 100755^6+63397^6=93338^6+82068^6+66612^6+62019^6+20469^6
     (MAVRIK and Greeri)
164. 100789^6+5599^6=98115^6+72793^6+42896^6+30786^6+26244^6
     (Liuqyn and Microcruncher*)
165. 101159^6+42058^6=97678^6+74046^6+59367^6+6760^6+6534^6
     (rroonnaalldd and Pwrguru)
166. 101183^6+77824^6=104102^6+49153^6+45024^6+19824^6+12390^6
     (Dietrich.li and tr0ach)
167. 103211^6+7033^6=97496^6+81795^6+60687^6+27086^6+22302^6
     (Jean-Charles Meyrignac and Laurent Evrard)
168. 104869^6+85444^6=103131^6+88634^6+53130^6+45342^6+10054^6
     (proteino.de and JagDoc)
169. 104989^6+46293^6=95487^6+90626^6+56865^6+38472^6+18396^6
     (proteino.de and Jeff17)
170. 106055^6+15685^6=103992^6+72723^6+40274^6+39066^6+36785^6
     (sqweeser and Jeff17)
171. 107035^6+52081^6=105847^6+65112^6+58758^6+9555^6+2482^6
     ([AF>HFR>RR] Jim PROFIT and Greeri)
172. 107726^6+84205^6=110358^6+59872^6+59437^6+48342^6+46410^6
     (Contrapunto Bestiale and Liuqyn)
173. 107876^6+78145^6=100737^6+92698^6+70632^6+7490^6+1428^6
     (Liuqyn and Laurent Evrard)
174. 108007^6+54189^6=104013^6+83352^6+47257^6+19278^6+13488^6
     (Clay and [AF>HFR>RR] Jim PROFIT)
175. 108211^6+895^6=107211^6+65366^6+45465^6+23240^6+5388^6
     (Laurent Evrard and [AF>HFR>RR] julien76100)
176. 108333^6+84905^6=106947^6+88935^6+11080^6+10122^6+4914^6
     ([P3D] Crashtest and Laurent Evrard)
177. 112741^6+59519^6=113106^6+35721^6+33722^6+27384^6+25475^6
     (Laurent Evrard and Jeff17)
178. 114109^6+6709^6=105306^6+97209^6+20734^6+15817^6+10920^6
     (Laurent Evrard and arctic light)
179. 114361^6+12374^6=102552^6+98658^6+65977^6+64036^6+24528^6
     (proteino.de and [AF>HFR>RR] julien76100)
180. 114658^6+77545^6=112184^6+88382^6+47292^6+44919^6+22578^6
     (Laurent Evrard and Clay)
181. 115715^6+102778^6=122102^6+78390^6+50883^6+50260^6+546^6
     (Friedrich and Gennady Stolyarov II)
182. 115819^6+72481^6=107516^6+100170^6+39036^6+24269^6+10647^6
     (Jeff17 and Laurent Evrard)
183. 118057^6+85428^6=118300^6+84126^6+27447^6+23772^6+408^6
     (Laurent Evrard and Fornax)
184. 118261^6+92465^6=114464^6+92452^6+84609^6+70242^6+15099^6
     (Gennady Stolyarov II and Mumps [MM])
185. 119314^6+55961^6=109886^6+98777^6+73662^6+63588^6+27822^6
     (Laurent Evrard and Renich)
186. 121544^6+6733^6=116622^6+85946^6+81816^6+41587^6+19260^6
     (Laurent Evrard and mikkovi)
187. 122361^6+98125^6=118458^6+106470^6+54628^6+37023^6+19707^6
     ([P3D] Crashtest and [AF>HFR>RR] julien76100)
188. 122797^6+3884^6=122550^6+57330^6+41274^6+29960^6+2317^6
     (Laurent Evrard and Jeff17)
189. 123227^6+16441^6=111708^6+106421^6+68746^6+12558^6+6825^6
     (ERBrouwer and Laurent Evrard)
190. 125407^6+26711^6=121728^6+88823^6+63999^6+59934^6+55930^6
     (Friedrich and http://jkfs.petrsu.ru/)
191. 125807^6+100545^6=127881^6+83112^6+75992^6+68376^6+20727^6
     (Mr. Hankey and Jeff17)
192. 127508^6+64461^6=127206^6+71379^6+23688^6+16542^6+4550^6
     (Killer4U and SEARCHER)
193. 128295^6+96391^6=124593^6+102466^6+80223^6+67704^6+25266^6
     (w42 and proteino.de)
194. 128353^6+34312^6=120225^6+106414^6+29328^6+26252^6+23562^6
     (HansNandi and Laurent Evrard)
195. 129593^6+84923^6=116025^6+111108^6+90811^6+75152^6+60708^6
     (zglloo and dimych)
196. 129650^6+79435^6=115590^6+114884^6+77574^6+67938^6+19957^6
     (aqua_mac and jjwhalen)
197. 129793^6+129098^6=134806^6+120204^6+78414^6+73658^6+15681^6
     (Laurent Evrard and toms83)
198. 130273^6+38718^6=120708^6+106113^6+80598^6+66210^6+47896^6
     (JagDoc and holly)
199. 130721^6+70948^6=123456^6+102132^6+83692^6+63490^6+56511^6
     (Antar and Daniil Zaytsev)
200. 131287^6+16880^6=124040^6+97335^6+92008^6+52692^6+16746^6
     (Laurent Evrard and XSmeagolX)
201. 131483^6+96151^6=129716^6+102921^6+40186^6+25221^6+12726^6
     (toms83 and Wanderer0815)
202. 133092^6+81845^6=130326^6+99162^6+44565^6+10738^6+1380^6
     (Odicin and ritterm)
203. 135299^6+44033^6=133749^6+86064^6+46820^6+22743^6+2912^6
     (Clay and Kristian_P)
204. 136790^6+61263^6=118846^6+116118^6+93807^6+92994^6+43638^6
     (Clay and Liuqyn)
205. 137059^6+105725^6=137634^6+101073^6+72519^6+48272^6+39374^6
     (Thorvin and Marc1)
206. 137695^6+9841^6=137487^6+57400^6+54116^6+29505^6+5034^6
     (Grandma061 and Laurent Evrard)
207. 138433^6+83493^6=130557^6+112791^6+78412^6+71442^6+11700^6
     (atom330 and LukaszST)
208. 138655^6+51721^6=120432^6+119294^6+101097^6+70504^6+32991^6
     (bachmann and [P3D] Crashtest)
209. 139295^6+49656^6=135786^6+99474^6+59388^6+58264^6+15393^6
     (firstomega and Wanderer0815)
210. 140450^6+73609^6=121512^6+118986^6+105266^6+86458^6+15561^6
     (Laurent Evrard and marcel kaiser)
211. 140915^6+4439^6=123375^6+121836^6+100128^6+53896^6+1925^6
     (UL1 and wn1hgb)
212. 141149^6+73727^6=132711^6+117000^6+58555^6+13650^6+8372^6
     ([AF>HFR>RR] julien76100 and Laurent Evrard)
213. 141195^6+48791^6=130786^6+118566^6+73227^6+17316^6+4137^6
     (RAZIELakaALIN and firstomega)
214. 141934^6+118033^6=132642^6+118440^6+106078^6+103761^6+25526^6
     (Laurent Evrard and [P3D] Crashtest)
215. 144663^6+96829^6=133203^6+123815^6+96270^6+41076^6+16380^6
     (Bent Vangli and King Lin)
216. 144937^6+38756^6=129946^6+127176^6+77832^6+41783^6+8820^6
     (JagDoc and refler)
217. 144986^6+29737^6=139671^6+104552^6+80622^6+73094^6+71370^6
     (Soul_keeper and Megacruncher of The Scottish Boinc Team)
218. 145073^6+3719^6=144313^6+77364^6+64911^6+12440^6+1848^6
     (Wanderer0815 and Conan)
219. 145249^6+141235^6=139622^6+123634^6+122919^6+119406^6+7317^6
     (L@MiR and serhjo)
220. 145387^6+39824^6=130200^6+120638^6+98568^6+91221^6+28994^6
     (toms83 and tom)
221. 147815^6+44635^6=132718^6+126147^6+96705^6+70824^6+2996^6
     (Laurent Evrard and Kristian_P)
222. 148655^6+32899^6=144221^6+106896^6+81438^6+46703^6+18396^6
     (Vato and Mumps [MM])
223. 150379^6+80382^6=148491^6+91602^6+85242^6+71890^6+33954^6
     (UL1 and [ESL] Jabberwoky)
224. 150440^6+81601^6=139846^6+126465^6+74592^6+71386^6+48060^6
     ([AF>HFR>RR] mirouf71 and guru2001-muc)
225. 152554^6+84143^6=151926^6+80080^6+78673^6+71736^6+54408^6
     (Jean-Charles Meyrignac and Toni)
226. 153775^6+77367^6=154077^6+57330^6+53018^6+23181^6+14760^6
     (PSL and proteino.de)
227. 154844^6+34235^6=148379^6+120260^6+66750^6+22176^6+15912^6
     (kwaa and Laurent Evrard)
228. 154871^6+135983^6=150255^6+128988^6+110208^6+110054^6+87191^6
     (Jeff17 and Mumps [MM])
229. 155521^6+105157^6=146958^6+132530^6+41272^6+37149^6+36279^6
     (john and Truth?)
230. 155789^6+119125^6=134313^6+132216^6+120905^6+118664^6+54180^6
     (dygon and bachmann)
231. 156825^6+36661^6=146874^6+125511^6+96075^6+71960^6+42132^6
     (Nasicus and wn1hgb)
232. 156895^6+72848^6=155526^6+85761^6+79674^6+72424^6+69790^6
     (Josef Andersson and Cocagne)
233. 157796^6+12107^6=143094^6+137328^6+71794^6+45045^6+2072^6
     (Soul_keeper and Fornax)
234. 157819^6+75763^6=156156^6+101815^6+49500^6+47815^6+7938^6
     (Iki and Mumps [MM])
235. 158075^6+98173^6=158249^6+87255^6+78774^6+69384^6+22054^6
     (mikkovi and toms83)
236. 158906^6+158509^6=159726^6+149100^6+127848^6+35306^6+26701^6
     (Buscho and Aflatoxin)
237. 159009^6+117545^6=141708^6+132378^6+132125^6+35022^6+2877^6
     (Infomat and Laurent Evrard)
238. 159373^6+109999^6=158536^6+113925^6+67431^6+30192^6+12122^6
     (Cazamarcianos and tom)
239. 159409^6+94135^6=152805^6+124670^6+90210^6+65730^6+28091^6
     (Cry and Coleslaw)
240. 160311^6+103199^6=156663^6+121254^6+72639^6+60340^6+53796^6
     (FreeLarry and vaughan)
241. 160331^6+33041^6=150122^6+132951^6+51363^6+27968^6+12012^6
     (kwaa and [SG-SPEG]atti)
242. 161011^6+3104^6=157395^6+109044^6+89316^6+56042^6+18116^6
     (L@MiR and Laurent Evrard)
243. 161998^6+109715^6=163422^6+95571^6+41990^6+36002^6+21000^6
     (DannyBoy and [FVG] mephyst)
244. 162159^6+20471^6=159411^6+99203^6+96264^6+52290^6+37764^6
     (GrafZahl and ????)
245. 163475^6+81565^6=155862^6+130909^6+44370^6+38381^6+5292^6
     (Mumps [MM] and DoctorNow)
246. 164116^6+158083^6=170030^6+134058^6+130186^6+82173^6+35634^6
     (Arcindo and [ESL Brigade] RinDorBroT)
247. 164117^6+106899^6=162396^6+114387^6+86124^6+58898^6+4725^6
     (bachmann and Jrachi)
248. 164712^6+106577^6=157794^6+133626^6+81615^6+45584^6+27090^6
     (DKD and Wassertropfen)
249. 165434^6+138177^6=157050^6+150465^6+96908^6+50820^6+42036^6
     (guru2001-muc and schensi)
250. 165533^6+76375^6=150173^6+142296^6+99708^6+51273^6+4774^6
     (yoyo[Friends] and Laurent Evrard)
251. 165806^6+29593^6=160615^6+123732^6+53076^6+7092^6+1414^6
     (L@MiR and Laurent Evrard)
252. 165909^6+88609^6=155120^6+139566^6+50769^6+11973^6+2604^6
     (Laurent Evrard and RushPill)
253. 166925^6+43514^6=162456^6+120589^6+72720^6+57070^6+9408^6
     (Jeff17 and Clay)
254. 168101^6+30904^6=151557^6+146926^6+83202^6+61562^6+27468^6
     (Mitchell and King Lin)
255. 168125^6+145417^6=173481^6+124712^6+99907^6+53508^6+13914^6
     (GLaDOS and Laurent Evrard)
256. 169871^6+41725^6=158277^6+142026^6+68180^6+38556^6+17495^6
     (dave_borg123 and SEARCHER)
257. 170117^6+94978^6=154686^6+145994^6+97371^6+93828^6+61544^6
     (Laurent Evrard and toms83)
258. 170117^6+144321^6=162972^6+138250^6+134979^6+106848^6+51849^6
     ([AF>France>Ouest>Normandie]The Stress Man (-: and [AF>Libristes] nico8313)
259. 170129^6+120725^6=168444^6+122219^6+93912^6+86352^6+63889^6
     (UL1 and Gennady Stolyarov II)
260. 170315^6+117800^6=154869^6+151536^6+102598^6+40992^6+37030^6
     (King Lin and Jean-Charles Meyrignac)
261. 171155^6+154453^6=173208^6+150555^6+57344^6+54020^6+42147^6
     (Laurent Evrard and Liuqyn)
262. 171209^6+83560^6=167164^6+113022^6+107122^6+69111^6+35364^6
     (Michelle and toms83)
263. 171236^6+44725^6=167892^6+118860^6+33936^6+20300^6+331^6
     (???? and SEARCHER)
264. 171286^6+56659^6=163730^6+124179^6+104774^6+99666^6+61152^6
     (Creatormaster and Laurent Evrard)
265. 173101^6+166379^6=162736^6+160875^6+146364^6+115507^6+14826^6
     (eric_lc and Laurent Evrard)
266. 173392^6+114573^6=171934^6+113820^6+106080^6+41811^6+7044^6
     ([AF>HFR>RR] ginie and Laurent Evrard)
267. 173823^6+167768^6=180402^6+152481^6+117180^6+79488^6+16912^6
     (MPG and Laurent Evrard)
268. 174745^6+74327^6=152574^6+141690^6+140329^6+81319^6+45864^6
     (Laurent Evrard and mackman123)
269. 175669^6+21473^6=167412^6+139505^6+32613^6+29524^6+17976^6
     (wn1hgb and Laurent Evrard)
270. 176023^6+98960^6=175959^6+99204^6+60562^6+33862^6+12012^6
     (bachmann and super-olli)
271. 176615^6+159991^6=186627^6+130170^6+38486^6+34965^6+23984^6
     (Kristian_P and http://urfak.petrsu.ru/)
272. 177901^6+55763^6=172272^6+124243^6+97552^6+96138^6+80013^6
     ([AF>WildWildWest] lamoule and Laurent Evrard)
273. 178249^6+60653^6=162988^6+146335^6+123012^6+67314^6+16149^6
     ([B^S] Ralfy and NGS~StugIII)
274. 179125^6+83065^6=161526^6+155865^6+98252^6+84294^6+39403^6
     (Cry and Coleslaw)
275. 179243^6+139045^6=179928^6+128632^6+103600^6+85605^6+81771^6
     (???? and RushPill)
276. 181389^6+132898^6=185544^6+82488^6+44520^6+35987^6+23724^6
     ([SG]marodeur6 and Jean-Charles Meyrignac)
277. 182435^6+16809^6=181056^6+108045^6+59836^6+39207^6+27630^6
     (Rumsey and Laurent Evrard)
278. 182857^6+27872^6=176982^6+123720^6+117166^6+88452^6+5257^6
     ([XTBA>XTC] Pousse Mousse and Clay)
279. 184103^6+16133^6=172095^6+148886^6+103971^6+91284^6+77980^6
     ([AF>Le_Pommier>MacBidouille.com] emeuley and wn1hgb)
280. 184409^6+16315^6=183134^6+104997^6+75975^6+63854^6+39060^6
     (guru2001-muc and whizbang)
281. 184799^6+47533^6=174282^6+143143^6+115563^6+97002^6+11732^6
     (Bold_Seeker and Reddogg)
282. 185254^6+16315^6=168084^6+153440^6+129892^6+45861^6+44280^6
     (mongfevned and Saxen)
283. 185362^6+29545^6=183834^6+103892^6+94276^6+43113^6+1152^6
     (Mr. Hankey and 7ri9991 [MM])
284. 186850^6+146159^6=191336^6+118608^6+77616^6+77574^6+51217^6
     (SETIKAH and Jeff17)
285. 186997^6+38987^6=173418^6+141057^6+140441^6+46920^6+19018^6
     (tom and Laurent Evrard)
286. 187184^6+13003^6=180888^6+141141^6+61136^6+53550^6+34412^6
     (Laurent Evrard and Matthew Hill)
287. 187239^6+34340^6=183582^6+118020^6+107815^6+90144^6+33894^6
     (King Lin and [ESL] Jabberwoky)
288. 187573^6+17085^6=183002^6+125454^6+106779^6+91728^6+50253^6
     (King Lin and shka)
289. 187769^6+77335^6=159494^6+155064^6+151683^6+98688^6+91189^6
     (Mumps [MM] and Laurent Evrard)
290. 188929^6+115081^6=182025^6+139839^6+123138^6+87934^6+3416^6
     ([AF>HFR>RR] Jim PROFIT and Laurent Evrard)
291. 189725^6+154134^6=189021^6+147294^6+117390^6+108210^6+38452^6
     ([b@h] tomcat and SEARCHER)
292. 190177^6+124559^6=181611^6+154966^6+102234^6+74669^6+756^6
     (Friedrich and Laurent Evrard)
293. 190907^6+6316^6=182910^6+135100^6+130200^6+46213^6+24594^6
     (HOME and Bold_Seeker)
294. 191431^6+17593^6=179039^6+157188^6+102648^6+52696^6+12015^6
     (Jonathan and Laurent Evrard)
295. 191759^6+142241^6=179606^6+167895^6+111258^6+68471^6+57486^6
     (glennpat and Dj Ninja)
296. 194947^6+109491^6=193218^6+118888^6+109662^6+47073^6+42777^6
     (Laurent Evrard and Jeff17)
297. 196233^6+178714^6=191310^6+165612^6+158327^6+127344^6+41076^6
     (SEARCHER and Wanderer0815)
298. 196367^6+150985^6=179438^6+156169^6+153468^6+129777^6+122934^6
     (Luke@ZMQLD and tony_51)
299. 196853^6+68291^6=188162^6+138618^6+137529^6+60809^6+8400^6
     (dygon and SM6GXQ Peter Lindquist)
300. 197257^6+26155^6=197106^6+76196^6+64591^6+35469^6+20370^6
     (Odicin and chinazrw)
301. 197285^6+137576^6=176832^6+174651^6+137046^6+68936^6+61292^6
     (whizbang and Enric Surroca)
302. 197489^6+142444^6=199840^6+118236^6+102991^6+58506^6+52962^6
     (Dietrich.li and Mr. Hankey)
303. 198545^6+21373^6=187482^6+147679^6+138612^6+84504^6+22477^6
     (Keck_Komputers and [AF>HFR>RR] alipse)
304. 198621^6+89809^6=167083^6+164892^6+164760^6+62457^6+36540^6
     (Jean-Charles Meyrignac and Laurent Evrard)
305. 198791^6+155153^6=185815^6+171807^6+139776^6+104812^6+27636^6
     (Akanee and Laurent Evrard)
306. 199096^6+127457^6=170244^6+167412^6+147071^6+147000^6+7462^6
     (toms83 and Laurent Evrard)
307. 199858^6+163483^6=193100^6+172452^6+127845^6+82082^6+27720^6
     (RushPill and ????)
308. 201983^6+12122^6=181932^6+174622^6+109998^6+107144^6+39585^6
     (Laurent Evrard and ttrue)
309. 202534^6+144017^6=205544^6+107569^6+90588^6+86898^6+35490^6
     (mrlumu and Laurent Evrard)
310. 202618^6+51257^6=181895^6+178938^6+71024^6+52122^6+52038^6
     (Jeff17 and Mr. Hankey)
311. 204649^6+42743^6=177384^6+162721^6+153720^6+148113^6+8428^6
     (Jeff17 and Mr. Hankey)
312. 205189^6+181399^6=214599^6+136250^6+127659^6+110782^6+50154^6
     (Creatormaster and Laurent Evrard)
313. 205201^6+35317^6=186277^6+177471^6+108560^6+25704^6+25302^6
     (http://jkfs.petrsu.ru/ and Jethro2000)
314. 206023^6+16063^6=200872^6+143598^6+108768^6+84105^6+27631^6
     (Liuqyn and Laurent Evrard)
315. 206226^6+47735^6=192192^6+163191^6+137946^6+94248^6+62344^6
     (Ananas and Administrator)
316. 206231^6+117353^6=204054^6+130242^6+108150^6+89437^6+84539^6
     (Mr. Hankey and Laurent Evrard)
317. 207829^6+118453^6=193896^6+175539^6+99043^6+38790^6+14308^6
     (Tex1954 and [AF>Libristes] e.clement)
318. 208549^6+90908^6=194390^6+166875^6+139230^6+22942^6+126^6
     (Laurent Evrard and Luke@ZMQLD)
319. 208763^6+97121^6=208488^6+105399^6+60957^6+60454^6+52220^6
     (Clay and Jeff17)
320. 209701^6+199073^6=211617^6+171780^6+168805^6+143246^6+11130^6
     (Truth? and Mr. Hankey)
321. 211145^6+4264^6=209815^6+121954^6+43512^6+11040^6+756^6
     (Bold_Seeker and Administrator)
322. 211733^6+73729^6=203007^6+157269^6+130806^6+70922^6+15260^6
     (Creatormaster and Administrator)
323. 213931^6+210085^6=216186^6+207207^6+91375^6+52920^6+26054^6
     (Fogle.SETI.USA [BlackOps] and Laurent Evrard)
324. 215071^6+181148^6=187551^6+186878^6+183078^6+148050^6+3326^6
     (Ananas and Odicin)
325. 215483^6+195649^6=229011^6+134680^6+130137^6+101122^6+59796^6
     (HSchmirPo and TBirdTheYuri)
326. 215891^6+20446^6=191181^6+187180^6+131286^6+116802^6+109844^6
     (Ananas and Fizban)
327. 215980^6+33193^6=213204^6+131142^6+115052^6+70889^6+60438^6
     (jepeeg and Daniil Zaytsev)
328. 215993^6+87373^6=214389^6+127608^6+88305^6+67438^6+1118^6
     (Bold_Seeker and Laurent Evrard)
329. 218369^6+76597^6=214690^6+145747^6+101961^6+36426^6+23058^6
     (Nutcracker and [AF>HFR>RR] julien76100)
330. 219505^6+21499^6=204666^6+181858^6+112935^6+69384^6+12355^6
     (Mr. Hankey and Mumps [MM])
331. 220610^6+159121^6=225302^6+92796^6+65328^6+2394^6+49^6
     ([AF>Le_Pommier>MacBidouille.com]Prof and Peter)
332. 220819^6+30665^6=206451^6+182672^6+98306^6+86121^6+58818^6
     (Major and Cry)
333. 222868^6+99537^6=217464^6+154511^6+126036^6+70056^6+52458^6
     (Ananas and [AF>HFR>RR] Jim PROFIT)
334. 222967^6+123911^6=219954^6+149394^6+109281^6+86879^6+29624^6
     (whizbang and Laurent Evrard)
335. 223015^6+166996^6=226380^6+134841^6+120586^6+100590^6+31058^6
     (Gemjunkie and toms83)
336. 223985^6+62597^6=211911^6+179634^6+105181^6+95970^6+59096^6
     (ekim7645 and Laurent Evrard)
337. 224339^6+105892^6=211092^6+179970^6+135336^6+81172^6+24371^6
     (TLC and serhjo)
338. 224396^6+163145^6=216132^6+172818^6+161304^6+82376^6+54551^6
     (Friedrich and [XTBA>XTC] Pousse Mousse)
339. 224533^6+100949^6=216282^6+158424^6+148421^6+83482^6+31515^6
     (Soul_keeper and Friedrich)
340. 225426^6+195919^6=212136^6+209538^6+149118^6+100177^6+4368^6
     (Ben and [AF>Libristes] Pascal)
341. 226018^6+118411^6=208580^6+185010^6+137109^6+123620^6+122658^6
     ([SG-2nd Wave]Ruputer and Laurent Evrard)
342. 226486^6+113297^6=202704^6+201768^6+74586^6+64988^6+20335^6
     (destroyfx and mikkovi)
343. 226705^6+93095^6=199119^6+198471^6+149212^6+99600^6+99148^6
     (HSchmirPo and changoo69)
344. 226801^6+221993^6=228221^6+204525^6+185822^6+70530^6+36960^6
     (Michaelis and Creatormaster)
345. 227515^6+184847^6=211644^6+204111^6+159397^6+1722^6+28^6
     (Tarkus_AN and BrunoROT)
346. 227875^6+177345^6=235165^6+105924^6+91302^6+18081^6+13272^6
     (proteino.de and [AF>HFR>RR] Jim PROFIT)
347. 228232^6+172379^6=230580^6+157157^6+112668^6+71106^6+59366^6
     ([SG] dingdong and wn1hgb)
348. 230857^6+105964^6=218526^6+186953^6+97566^6+74502^6+73990^6
     (Heijo and Laurent Evrard)
349. 231680^6+166387^6=218484^6+183750^6+173523^6+104428^6+37360^6
     (Mitchell and Conan)
350. 231918^6+130745^6=213927^6+200382^6+43314^6+25040^6+16800^6
     (Wanderer0815 and firstomega)
351. 233577^6+1543^6=220972^6+183894^6+138999^6+67431^6+2646^6
     (firstomega and wn1hgb)
352. 234083^6+65375^6=217924^6+177408^6+172458^6+35469^6+21743^6
     (Norman_RKN and Fornax)
353. 234159^6+103952^6=234432^6+68474^6+31695^6+31038^6+26754^6
     (Mark Roman Schellberg and Rebecca)
354. 234299^6+129697^6=213264^6+200979^6+131257^6+130830^6+65632^6
     (Laurent Evrard and Yochanon)
355. 234347^6+151401^6=226085^6+167769^6+166218^6+95130^6+42420^6
     (Major and Conan)
356. 235123^6+86285^6=232323^6+144816^6+119462^6+31364^6+2877^6
     (Mumps [MM] and Phitar)
357. 235352^6+198157^6=241276^6+165166^6+145299^6+123186^6+27462^6
     (Serge Vannieuwenborgh and DoctorNow)
358. 235936^6+178727^6=241974^6+121017^6+103026^6+45920^6+45178^6
     (jjwhalen and whizbang)
359. 236363^6+24349^6=213423^6+206800^6+108564^6+49266^6+7217^6
     (Laurent Evrard and Kristian_P)
360. 237425^6+85667^6=233454^6+158280^6+103768^6+90685^6+68943^6
     (ongekend41 and gow)
361. 237667^6+80215^6=217665^6+204834^6+70171^6+66066^6+65516^6
     (proteino.de and johnnymcknowledge)
362. 237836^6+144025^6=201292^6+193860^6+185703^6+175588^6+18228^6
     (Deltanuss and Jeff17)
363. 238768^6+180709^6=228396^6+204038^6+122430^6+117598^6+8379^6
     (Trog Dog and Ananas)
364. 239698^6+186647^6=211575^6+204132^6+194222^6+147546^6+134408^6
     (LaVey and JagDoc)
365. 239707^6+163781^6=242921^6+106740^6+106449^6+86508^6+74102^6
     (SEARCHER and NGS~StugIII)
366. 239887^6+66554^6=236224^6+158865^6+96642^6+36792^6+20366^6
     ([B^S] Spydermb and [AF>Le_Pommier>MacADSL.com]Bertrand)
367. 239942^6+53987^6=239811^6+89646^6+67522^6+56644^6+36036^6
     (RushPill and tom)
368. 240067^6+239173^6=241322^6+237650^6+97965^6+64428^6+22635^6
     (lcdad78039[SETI.USA] and Gennady Stolyarov II)
369. 240187^6+19366^6=230706^6+170235^6+151462^6+129626^6+62412^6
     (???? and Wanderer0815)
370. 240197^6+22891^6=238210^6+129414^6+128772^6+64601^6+43953^6
     (Jeff17 and dusti)
371. 240797^6+127153^6=239896^6+142477^6+75978^6+40908^6+735^6
     (ExtraTerrestrial Apes and Clay)
372. 242903^6+8333^6=240261^6+146720^6+117393^6+82972^6+70812^6
     (Miri37 and LordNord)
373. 243062^6+144729^6=235648^6+173817^6+159678^6+48342^6+16332^6
     (JPLiz and Laurent Evrard)
374. 243685^6+180302^6=240226^6+169092^6+161130^6+144645^6+107408^6
     (Creatormaster and RushPill)
375. 243781^6+4509^6=238308^6+170532^6+111588^6+76979^6+16107^6
     (Laurent Evrard and [AF>Le_Pommier>MacBidouille.com] emeuley)
376. 247261^6+232214^6=249704^6+212492^6+192477^6+38682^6+9408^6
     (sleeplift and cze_siek)
377. 247449^6+98411^6=246708^6+129549^6+80532^6+36543^6+3892^6
     (GrafZahl and [SG]KidDoesCrunch)
\end{verbatim}

\begin{footnotesize}

\end{footnotesize} 
\end{document}